\documentclass{amsart}
\usepackage[all]{xy}
\newtheorem{innercustomtheorem}{Theorem}
\newenvironment{customtheorem}[1]
  {\renewcommand\theinnercustomtheorem{#1}\innercustomtheorem}
  {\endinnercustomtheorem}

\def\P{\mathbb{P}}

\def\F{\mathcal{F}}

\def\H{\mathcal{H}}

\def\I{\mathcal{I}}

\def\E{\mathcal{E}}

\def\dim{\mathrm{dim}}

\def\separation{\medskip}
\def\acc{\`}
\theoremstyle{plain}
\newtheorem{theorem}{Theorem}[section]
\newtheorem{lemma}[theorem]{Lemma}
\newtheorem{prop}[theorem]{Proposition}

\newtheorem{corollary}[theorem]{Corollary}
\newtheorem{remarks}[theorem]{Remarks}
\newtheorem{remark}[theorem]{Remark}

\newtheorem{example}[theorem]{Example}

\begin{document}

\title{The local cohomology of the jacobian ring}
\date{}
\author{Edoardo Sernesi}
 \maketitle
\begin{abstract}
 We study the 0-th local cohomology module $H^0_{\mathbf{m}}(R(f))$ of 
the jacobian ring $R(f)$ of a \emph{singular}  reduced complex projective hypersurface $X$, by relating it to the sheaf of logarithmic vector fields along $X$.
We investigate the analogies between $H^0_{\mathbf{m}}(R(f))$ and   the well known properties of the jacobian ring of a nonsingular hypersurface. 
 In particular we   study     self-duality, Hodge theoretic  and Torelli type questions for $H^0_{\mathbf{m}}(R(f))$.  

\end{abstract}

\section{Introduction}

 In this paper we focus on   the relation existing
 between   a (singular) projective hypersurface and the $0$-th local
cohomology of its jacobian ring. Most of the  results we will present are   well known to the experts, 
and perhaps the only novelty here is a unifying approach obtained by relating the local cohomology 
to the sheaf of logarithmic vector fields along $X$.
We will take the opportunity to introduce  what seem to us some interesting open
problems on this subject.

 Consider the polynomial ring
$P=\mathbb{C}[X_0,\dots,X_r]$ in $r+1$ variables, $r \ge 2$, with
coefficients in $\mathbb{C}$. Given a reduced polynomial $f \in P$
homogeneous of degree $d$ let $X := V(f) \subset \P^r$ be the
hypersurface defined by $f$. The \emph{jacobian ring} of $f$ is
defined as
\[
R=R(f) := P/J(f)
\]
where
\[
J(f) := \left({\partial f \over \partial X_0}, \dots, {\partial f
\over \partial X_r} \right)
\]
is the  \emph{gradient ideal} of $f$.

 If $X$ is nonsingular
then   $J(f)$ is generated by a regular sequence, and $R(f)$ is a
Gorenstein artinian ring with socle in degree $\sigma:=(r+1)(d-2)$.
It carries information on the geometry of $X$ and on its period
map. This classical case has been studied by Griffiths  
and his school. In   \cite{pG69} Griffiths has
shown the relation existing between the jacobian ring of a
nonsingular projective hypersurface and the Hodge decomposition of
its  primitive cohomology in middle dimension, and studied the
relation of $R(f)$ with the period map
    (see also \cite{cV03} for details and \cite{dC87} for
a  survey).

Assume now that   $X \subset \P^r$ is  singular, but reduced. In
this case the jacobian ring is not of finite length, in particular 
it is not artinian Gorenstein any more. It contains information  
on the structure of the singularities and on the global geometry
of $X$. This situation has been studied extensively,
both from the point of view of singularity theory (see e.g.
\cite{mG80,jS85,lT86,jW87,WY07}) and in relation with the (mixed) Hodge
theory of   $U:=\P^r\setminus X$ (see \cite{pD71,CD94,DD90,aD92,DS11}).
  Our main purpose is  to indicate a method to
distinguish the global information contained in $R(f)$ from the
local one coming from the nature of the singularities.

  Our starting point is the observation that, if $X$ is nonsingular,
   we have a canonical
  identification of $P$-modules
  \[
  R(f) = \bigoplus_{j\in \mathbb{Z}}H^1(T\langle X\rangle(j-d))
  \]
  where $T\langle X\rangle$ is the subsheaf of $T_{\P^r}$ of
  logarithmic vector fields   along $X$. If $X$ is
  singular this identification does not hold,
  but the $P$-module on the right hand side is the $0$-th local
  cohomology of $R(f)$.    We will see that this object contains
    relevant global informations about $X$.

To any finite type graded $P$-module one can associate
a coherent sheaf $M^{\sim}$ on $\P^r$ and there is a well-known
exact sequence involving the local cohomology graded modules (see \cite{EGA}, Prop. 2.1.5):
\begin{equation}\label{E:int1}
    \xymatrix{
    0\ar[r]&H^0_{\mathbf{m}}(M)\ar[r]&M\ar[r]&
    H^0_*(M^\sim)\ar[r]&H^1_{\mathbf{m}}(M)\ar[r]&0
    }
\end{equation}
where we used the notation $H^i_*(\F)=\bigoplus_{\nu\in
\mathbb{Z}}H^i(\F(\nu))$ for a coherent sheaf $\F$. In case
$M=R(f)$ with $X$ singular both $H^i_{\mathbf{m}}(R(f))$ are
finite length modules that carry interesting information about the
hypersurface $X$. In particular, $H^0_{\mathbf{m}}(R(f))$ contains
global information about $X$, while $H^1_{\mathbf{m}}(R(f))$ is
related with the singularities of $X$. We want to collect evidence
supporting the following principle:

\medskip

\emph{Most properties of the jacobian ring $R(f)$ in the
nonsingular case are transferred to  the local cohomology module
$H^0_{\mathbf{m}}(R(f))$  if $X$ is a singular hypersurface.}

\medskip

In particular one expects the following   in some
generality:

\begin{itemize}
    \item[(a)] \emph{Self-duality}, extending the analogous
    property of Artinian Gorenstein algebras.
    \item[(b)] \emph{Existence of  a connection with moduli of $X$,} in particular
    with first order locally trivial deformations of $X$.
    \item[(c)]  \emph{Existence of a relation with the Hodge decomposition} of the middle
    dimension primitive cohomology a nonsingular model of $X$.
		
		\item[(d)] \emph{Torelli type results}, stating the
possibility of reconstructing $X$ from $H^0_{\mathbf{m}}(R(f))$,
under some hypothesis. 
\end{itemize}

Question (a) has already attracted the attention of several authors and some results are  known. 
One is led naturally to consider more generally
the $0$-th local cohomology of algebras of the form $P/I$ where
$I=(f_0,\dots,f_r)$ is an ideal generated by $r+1$ homogeneous
polynomials, of degrees $d_0, \dots, d_r$. The following result is a special case of \cite{vSW11}, Theorem 4.7:

\begin{customtheorem}{\ref{C:spseq}}
 Assume that $\mathrm{dim}[\mathrm{Proj}(P/I)]=0$. Then
 there is a natural isomorphism:
 \[
 H^0_{\mathbf{m}}(P/I)\cong [H^0_{\mathbf{m}}(P/I)(\sigma)]^\vee
 \]
 where $\sigma = \sum_{i=0}^r d_i -r-1$. In particular we have
 natural isomorphisms
 \[
 H^0_{\mathbf{m}}(P/I)_k \cong H^0_{\mathbf{m}}(P/I)_{\sigma-k}^\vee, \quad 0 \le k \le \sigma
 \]
\end{customtheorem}

We include an independent
proof  of Theorem \ref{C:spseq}, more related with our point of view,  which uses a spectral sequence argument and is an adaptation of
the standard proof of Macaulay's Theorem (see e.g. \cite{cV03}).  
I am also aware of work in progress of H. Hassanzadeh and A. Simis about extensions of
Theorem \ref{C:spseq} to a local algebra situation.  
Taking $I=J(f)$, as a special case we obtain:

\begin{customtheorem}{\ref{T:selfdu}} 
    Assume that the hypersurface $X$ has only isolated
    singularities.  Then:
    \[
H^0_{\mathbf{m}}(R(f))_k \cong
H^0_{\mathbf{m}}(R(f))_{\sigma-k}^\vee, \quad 0 \le k \le \sigma
 \]
where $\sigma = (r+1)(d-2)$.
\end{customtheorem}

This is a   generalization of Macaulay's Theorem, that
states the self-duality of $R(f)$ in the case $X$ nonsingular. 
The theorem, in an equivalent form, appeared already in \cite{DS12}, Theorem 1. 
A similar result for hypersurfaces with isolated quasi-homogeneous singularities is proved in \cite{EM13}.
We also refer the reader to the recent preprint \cite{EU13} where all these duality 
results are reconsidered and further generalized. For recent related work see \cite{gS131,gS132}.

As mentioned before, we interpret $H^0_{\mathbf{m}}(R(f))$ by
means of the sheaf $T\langle X\rangle$, also denoted by
 $\mathrm{Der}(-\mathrm{log} X)$, associated to
any hypersurface $X$ in a smooth variety $M$ (see \S \ref{S:log}
where we recall its definition). Precisely we show that there is
an identification:
\begin{equation}\label{E:int2}
H^0_{\mathbf{m}}(R(f))=H^1_*(T\langle X\rangle(-d))
\end{equation}
(Proposition \ref{P:log1}). In particular:
\begin{equation}\label{E:int3}
H^0_{\mathbf{m}}(R(f))_d= H^1(T\langle X\rangle)
\end{equation}
The right hand side is the space of first-order locally trivial
deformations of $X$ in $\P^r$  (see \cite{eS06}, \S 3.4.4).
Therefore (\ref{E:int3}) generalizes what happens in the
nonsingular case, when we have the identification of $R(f)_d$ with
the space of first order deformations of $X$ in $\P^r$ modulo
projective automorphisms \cite{dC87}. Thus (\ref{E:int3}) gives an
answer to (b).

 In passing
note that Theorem \ref{T:selfdu} and (\ref{E:int2}) together imply the
self-duality of $H^1_*(T\langle X\rangle(-d))$ in the case when
$X$ has isolated singularities. This fact is quite straightforward
when $r=2$  
but it is not so  when $r \ge 3$, since $T\langle
X\rangle$ is not even locally free.

\separation

As of Question (c), one
 expects that there exists a   relation between the local cohomology of $R(f)$ 
 and the Hodge decomposition of the middle primitive cohomology of a nonsingular model
 $X'$ of $X$. We collect some evidence that this relation exists at least for strictly 
normal crossing hypersurfaces. In particular we show that for such hypersurfaces one has an isomorphism
\[
H^0_{\mathbf{m}}(R(f))_{d-r-1} \cong \bigoplus_{i=1}^s H^0(X_i,\Omega^{r-1}_{X_i})
\]
where $X_1, \dots, X_s$ are the irreducible components of $X$ (see Theorem \ref{T:ncd2} for a precise statement). This result and duality
imply a result completely analogous to Griffiths' for strictly normal crossing plane curves 
 (Corollary \ref{C:ncd1}).  We also prove a result for surfaces in $\P^3$ indicating that   
   the local cohomology contains  information on how the various components intersect (Theorem
\ref{T:ncd3}).

Question (d) is related to interesting issues that have been widely considered 
in the case of arrangements of hyperplanes and of hypersurfaces, but from a 
different point of view. Several authors have  investigated the problem of
reconstructing certain arrangements of hyperplanes and of hypersurfaces  from their sheaf
of logarithmic differentials (see \cite{eA13,DK93,FMV,UY09,jV00}). Our Question (d) is quite different, at least when $r \ge 3$, while it is 
 essentially equivalent to it when $r=2$.  We discuss the problem and we give a few examples.     

In the paper we also consider the question of freeness of the sheaf $T\langle X \rangle$, which is a special case of  
the condition  $H^0_{\mathbf{m}}(R(f))=0$.  We overview some of the known results     in
the case $r=2$.

In detail the paper is organized as follows. \S \ref{S:log}  is devoted to
 the relation between local cohomology of the jacobian ring
of $X$ and the sheaf $T\langle X\rangle$.  
In \S \ref{S:local} we    
consider the self duality properties.   \S \ref{S:logar} is devoted to generalities on sheaves of logarithmic differentials and \S \ref{S:hodge}  
 to the Hodge theoretic properties of the local cohomology.
  In the next   \S \ref{S:torelli} we discuss the Torelli problem (d) above, and its relations with 
related reconstruction problems. \S \ref{S:free} treats the freeness of $T\langle X \rangle$.

\textbf{Acknowledgements.}    I am grateful to H. Hassanzadeh
and A. Simis for useful remarks concerning Theorem \ref{C:spseq}, to D. Faenzi for his
help with Example (\ref{Ex:simis}),
    to E. Arbarello, F. Catanese,  A. Lopez and J. Valles for
helpful conversations.  
All the examples have been computed using Macaulay2 \cite{GS94}.

  After posting the first version 
of this paper I became aware of references \cite{DS12} and \cite{vSW11}. 
I am thankful to D. Van Straten,  and M. Saito for calling my attention on them and for some helpful remarks.
Finally it is a pleasure to thank A. Dimca for his correspondence and for bringing 
Example \ref{Ex:redsurf} to my attention. 

I am a member of
INDAM-GNSAGA.  This research has been supported by the project
MIUR-PRIN 2010/11 \emph{Geometria delle variet\`a algebriche.}

%%%%%%%%%%%%%%%%%%%%%%%%%%%%%%%%%%%%%%%%%%%%%%%%%%%%%%%%%%%%%%%%

\section{Logarithmic derivations and local cohomology}\label{S:log}

We will adopt the following standard notation and terminology. Consider the graded polynomial ring
$P=\bigoplus_{k\ge 0} P_k = \mathbb{C}[X_0,\dots,X_r]$, in $r+1$ variables, $r \ge 2$, with
coefficients in $\mathbb{C}$, and denote by $\mathbf{m}=\bigoplus_{k\ge 1}P_k$   its irrelevant maximal ideal. 
  A graded $P$-module $M = \bigoplus_k M_k$ is \emph{$TF$-finite} if $M_{\ge k_0}:=\bigoplus_{k\ge k_0} M_k$
	is of finite type for some $k_0$. If $M$ is $TF$-finite   we let 
\[
M^\vee = \bigoplus_k (M^\vee)_k = \bigoplus_k M_{-k}^\vee = \bigoplus_k \mathrm{Hom}_{\mathbb{C}}(M_{-k},\mathbb{C})
\]   
For any coherent sheaf $\mathcal{F}$ on $\mathbb{P}^r$
and $0\le i \le r$ we let
\[
H^i_*(\mathcal{F}) = \bigoplus_{k\in\mathbb{Z}} H^i(\mathbb{P}^r,\mathcal{F}(k))
\]
which is a  graded $P$-module.

Consider  a reduced polynomial $f \in P$
homogeneous of degree $d$.  Let $X := V(f) \subset \P^r$ be the
hypersurface defined by $f$ and let  
 \[
R(f) := P/J(f)
\]
be the \emph{jacobian ring}  of $f$ (or of $X$) where
\[
J(f) := \left({\partial f \over \partial X_0}, \dots, {\partial f
\over \partial X_r} \right)
\]
is the  \emph{gradient ideal} of $f$.  The scheme
$\mathrm{Proj}(R(f))$ is called the \emph{jacobian scheme}  of $f$, or the \emph{singular scheme}    of $X$ (see \cite{pA95}), 
and also denoted by  $\mathrm{Sing}(X)$.  
 We denote by $\mathcal{J}_f= J(f)^\sim \subset \mathcal{O}_{\P^r}$ the ideal sheaf associated to $J(f)$, and by
	\[
	\mathcal{J}_{f/X} = \mathcal{J}_f/\mathcal{I}_X \subset \mathcal{O}_X
	\]
	its image in $\mathcal{O}_X$. Then $\mathcal{J}_{f/X}$ is called \emph{the jacobian ideal sheaf of $X$}. 
	Note that  $\mathcal{O}_X/\mathcal{J}_{f/X}= \mathcal{O}_{\mathrm{Sing}(X)}=T^1_X(-d)$ where $T^1_X$ is  \emph{the first cotangent sheaf of $X$,} .
	
	A more useful description of  the jacobian ring is the following. Consider the diagram of sheaf homomorphisms:
	\begin{equation}\label{E:log0}
	\xymatrix{
	  0\ar[r]&\mathcal{K} \ar[r]\ar@{=}[d] &\mathcal{O}_{\P^r}(-d+1)^{r+1}\ar@{=}[d] \ar[r]^-{\partial f}& \mathcal{O}_{\P^r} \ar[r] &  T^1_X(-d) \ar[r] & 0 \\
		0\ar[r]&\mathcal{K} \ar[r]  &\mathcal{O}_{\P^r}(-d+1)^{r+1}  \ar[r]^-{\partial f}& \mathcal{J}_f \ar[r]\ar@{^(->}[u] & 0}
	\end{equation}
	where $\partial f$ is defined by the partials of $f$, and $\mathcal{K}=\mathrm{ker}(\partial f)$.  It induces  
	\[
	\xymatrix{
	H^0_*(\mathcal{O}_{\P^r}(-d+1))^{r+1} \ar[r]\ar@{=}[d] & P \ar[r] & R(f) \ar[r] &0 \\
	H^0_*(\mathcal{O}_{\P^r}(-d+1))^{r+1} \ar[r] & J(f)^{sat} \ar[r]\ar@{^(->}[u] & H^1_*(\mathcal{K})\ar@{^(->}[u] \ar[r] & 0
	}
	\]
	where $J(f)^{sat}$ is the saturation of $J(f)$. The following are clearly equivalent conditions:
	\begin{itemize}
		\item[(a)] $X$ is nonsingular.
		
		\item[(b)] $T^1_X=0$.
		
		\item[(c)] $R(f)$  has finite length.
		
		\item[(d)]  $R(f) = H^1_*(\mathcal{K})$.
	\end{itemize}
	
	When they are not satisfied then $H^1_*(\mathcal{K})$ is just a submodule of finite length of $R(f)$ 
	and we have an identification:
	\begin{equation}\label{E:log1}
	H^1_*(\mathcal{K})  = \frac{J(f)^{sat}}{J(f)} = H^0_{\mathbf{m}}(R(f))
	\end{equation}
	where $H^0_{\mathbf{m}}(M)$ denotes the \emph{$0$-th local cohomology} of a graded $P$-module  with respect to 
	$\mathbf{m}$.

 We also have the exact sequence:
	\begin{equation}\label{E:log2}
	\xymatrix{
	0 \ar[r]& T\langle X \rangle \ar[r] & T_{\P^r} \ar[r]^-\eta& \mathcal{J}_{f/X}(d) \ar[r] & 0}
	\end{equation}
	where $T\langle X\rangle:= \mathrm{ker}(\eta)$ is the
\emph{sheaf of logarithmic vector fields along $X$} and $\eta$ is defined as:
\[
\eta\left(\sum_i A_i(X)\frac{\partial}{\partial X_i}\right) = \sum_i A_i {\partial f \over \partial X_i}
\]
 the sheaf $T\langle X \rangle$ is
also denoted by $\mathrm{Der}(-\mathrm{log} X)$ in the literature
\cite{kS80}.
	We then have the following commutative diagram with exact rows and columns:
	\[
	\xymatrix{&&0&0 \\
	0\ar[r]& T\langle X \rangle \ar[r] & T_{\P^r}\ar[u] \ar[r]^-\eta& \mathcal{J}_{f/X}(d)\ar[u] \ar[r] & 0 \\
	0\ar[r]&\mathcal{K}(d)\ar[u]^-{\cong}\ar[r]&\mathcal{O}_{\P^r}(1)^{r+1} \ar[u]\ar[r]^-{\partial f} & \mathcal{J}_f(d) \ar[r]\ar[u] &0 \\
	&&\mathcal{O}_{\P^r}\ar[u]\ar[r]^-f& \mathcal{I}_X(d)\ar[u]\ar[r]&0 \\
	&&0\ar[u]&0\ar[u]}
\]
where the middle vertical is the Euler sequence. From this diagram we deduce the isomorphisms:
\begin{align}
T\langle X \rangle &\cong \mathcal{K}(d)\label{E:log3}\\
H^1_*(\mathcal{J}_{f/X})&\cong H^1_*(\mathcal{J}_f)\ (=H^2_*(\mathcal{K}) \ \text{if $r\ge 3$})\label{E:log4}
\end{align}
Now we can prove the following:

\begin{prop}\label{P:log1}
In the above situation we have a canonical isomorphism:
\begin{align}
H^0_{\mathbf{m}}(R(f)) &\cong H^1_*(T\langle X \rangle(-d))\label{E:log5}  \end{align}
In particular \[R(f) \cong H^1_*(T\langle X \rangle(-d))\] if $X$ is nonsingular.
\end{prop}

\proof It  follow directly from (\ref{E:log1}) and (\ref{E:log3}). The last assertion is obvious. \qed

\begin{corollary}\label{C:log2}
The vector space  $H^0_{\mathbf{m}}(R(f))_d$ is naturally identified with the space of
first order locally trivial deformations of $X$ in $\P^r$ modulo the action of $\mathrm{PGL}(r+1)$.
\end{corollary}

\proof The proposition identifies $H^0_{\mathbf{m}}(R(f))_d$ with  $H^1(T\langle X \rangle)$ which is the space of first order 
locally trivial deformations of the inclusion $X \subset \P^r$ (see \cite{eS06}, \S 3.4.4  p. 176).
\qed

\begin{remarks}\rm\label{R:log1}
(i) It is easy to compute  that for $X \subset \P^2$ the \emph{Chern classes of} $T\langle X\rangle(k)$ are:
\[
c_1(T\langle X\rangle(k)) = 3-d+2k, \quad c_2(T\langle
 X\rangle(k)) = d^2-(3+k)d+3+3k+k^2- t^1_X
\]
where       $t^1_X = h^0(T^1_X) = h^0(\mathcal{O}_{\mathrm{Sing}(X)})$.  Moreover:
\[
-\chi(T\langle X\rangle) =  \frac{1}{2}d(d+3)-t^1_X-8 
\]
which is the expected dimension of the family of locally trivial deformation of $X$ modulo $\mathrm{PGL}(3)$. 
This is explained by the fact that $T\langle X\rangle$ is the sheaf controlling the locally trivial deformation theory of $X$ in $\P^2$ 
(see \cite{eS06}).

\separation

(ii) If $X$ is a normal crossing arrangement of $d \ge r+2$ hyperplanes then $T\langle X\rangle$ is the dual of a 
	Steiner bundle \cite{DK93}, in particular it is locally free, and these bundles are known to be stable \cite{BS92}.  In the special case $d=r+2$ we have
	$T\langle X\rangle= \Omega(1)$. If $1 \le d \le r+1$ then
	\[
	T\langle X\rangle = \mathcal{O}_{\P^r}^{d-1} \bigoplus \mathcal{O}_{\P^r}(1)^{r+1-d}
	\]
	and these bundles are not stable.

(iii) If $X\subset \P^2$ is nonsingular then $T\langle X\rangle$ is stable (\cite{UY07}, Lemma 3).  

\end{remarks}

In the case of plane curves  we have more generally:

\begin{prop}\label{P:stable}
Let $X \subset \P^2$ be of degree $d\ge 4$.  Then $T\langle X\rangle$ is stable if and only if $(f_0,f_1,f_2)$, where
$f_i={\partial f \over \partial X_i}$, has no syzygies of degree $\left[(d-1)/2\right]$.
In particular $T\langle X\rangle$ is stable if $X$ is nonsingular.
\end{prop}

\proof Twist $T\langle X\rangle$ by $k=\left[(d-3)/2\right]$.  Then $c_1(T\langle X\rangle(k))=0,-1$ according to whether $d$ is odd or even, and 
  $T\langle X\rangle$ is stable if and only if $H^0(T\langle X\rangle(k))=0$  (\cite{OSS}, Lemma 1.2.5 p. 165).  The 
exact sequence
\[
\xymatrix{
0 \ar[r] & T\langle X\rangle(k) \ar[r] & \mathcal{O}_{\P^2}(k+1) \ar[r]^-{(f_0,f_1,f_2)} & \mathcal{J}_f(d+k) \ar[r] & 0}
\]
identifies $H^0(T\langle X\rangle(k))$ with the space of syzygies of $(f_0,f_1,f_2)$
of degree $k+1 = \left[(d-1)/2\right]$.  

In the nonsingular case $(f_0,f_1,f_2)$ has no syzygies of degree less than $d-1$ because they form a regular sequence. \qed

\begin{example}\label{Ex:unstable}\rm
Let $f= X_1^\alpha X_0^{d-\alpha} - X_2^d$,  with $2 \le \alpha < d$, and $d \ge 4$.  Then $T\langle X\rangle$ is not stable because $(f_0,f_1,f_2)$ has the 
linear syzygy $(\alpha X_0,-(d-\alpha)X_1,0)$.

Additional interesting informations concerning the syzygies of $(f_0,f_1,f_2)$ for a singular plane curve are in \cite{aD13}.
\end{example}

%%%%%%%%%%%%%%%%%%%%%%%%%%%%%%%%%%%%%%%%%%%%%%%%%%%%%%%%%%%%%%%%%%%%%%%%%

\section{Self-duality of the local cohomology}\label{S:local}

In this section we will consider a situation slightly more general than before. Let
\[
I = (\mathbf{f}) = (f_0, \dots, f_s) \subset P
\]
be a proper homogeneous ideal, whose generators have degrees $d_0,\dots,d_s$ respectively, and let $R = P/I$.
Denote by $Y = \mathrm{Proj}(R)$ and by $\mathcal{I} = I^\sim \subset \mathcal{O}_{\mathbb{P}^r}$.  We have an exact sequence:
\begin{equation}\label{E:locohom1}
\xymatrix{
 0 \ar[r] & \mathcal{K} \ar[r] & \bigoplus_{j=0,\dots,s} \mathcal{O}_{\mathbb{P}^r}(-d_j) \ar[r]^-{\mathbf{f}} & \mathcal{I} \ar[r] &0
}
\end{equation}
where $\mathcal{K}:= \mathrm{ker}(\mathbf{f})$. The \emph{$0$-th and $1$-st local cohomology modules of $R$} (with respect to $\mathbf{m}$)
are defined respectively as:
\begin{align}
H^0_{\mathbf{m}}(R)&:=H^1_*(\mathcal{K}) \notag \\
H^1_{\mathbf{m}}(R)&:=H^1_*(\mathcal{I})\ (= H^2_*(\mathcal{K}) \ \ \text{if $r\ge 3$})\notag
\end{align}
They are graded $P$-modules of finite length. In case $\mathbf{m}^k \subset I$  for some $k > 0$, i.e. $Y = \emptyset$,
 we have
\[
H^0_{\mathbf{m}}(R) = R, \ \ H^1_{\mathbf{m}}(R) = (0)
\]
  There is a standard exact sequence:
\begin{equation}\label{E:locohom2}
    \xymatrix{
    0\ar[r]& H^0_{\mathbf{m}}(R)\ar[r]& R \ar[r] &
    H^0_*(\mathcal{O}_Y)\ar[r] & H^1_{\mathbf{m}}(R)\ar[r]&0 }
\end{equation}
\emph{Assume now that $s=r$.}  Denote by     
\[
 \E := \bigoplus_{j=0,\dots,r} \mathcal{O}_{\P^r}(-d_j)
 \]
and let 
\[
\sigma := \sum_j (d_j-1) = \sum_j d_j -r-1
\]
 Consider the Koszul complex:
 \[
 \E^\bullet:\xymatrix{
 0\ar[r]&
 \E_{-r-1}\ar[r]&\E_{-r}\ar[r]&\cdots\ar[r]&\E_{-1}\ar[r]^-{\mathbf{f}}&\E_0\ar[r]&0
 }
 \]
 where $\E_{-p} =\bigwedge^p \E$.  For every $k \in \mathbb{Z}$ we can consider the twist
 $\E^\bullet(k)$ 
 and   the two corresponding spectral sequences of hypercohomology. Taking direct sums over all $k$ we can
 collect them in the following two   spectral
 sequences of hypercohomology:
 \[
  \begin{array}{l}
   A_1^{pq}= H^q_*(\E_p) \\
   B_2^{pq}= H^p_*(\H^q(\E^\bullet)) \\
 \end{array}
\]
where $\H^q(\E^\bullet)$ is the $q$-th cohomology sheaf of
$\E^\bullet$. In particular $\H^0(\E^\bullet)=\mathcal{O}_Y$.  In the
$A$-spectral sequence we have in particular:
\begin{equation}\label{E00}
A_2^{00}= \cdots =
A_{r+1}^{00}=\mathrm{coker}[H^0_*(\E_{-1})\longrightarrow
H^0_*(\E_0)] =R
\end{equation}
\begin{equation}\label{Err}
A_2^{-r-1r}=\cdots = A_{r+1}^{-r-1r}=
\mathrm{ker}[H^r_*(\E_{-r-1})\longrightarrow H^r_*(\E_{-r})] =
[R(\sigma)]^\vee
\end{equation}
and
\[ \mathbf{d}_{r+1}:[R(\sigma)]^\vee=A_{r+1}^{-r-1r}\longrightarrow A_{r+1}^{00}=R
\]
 We denote
by $\mathbf{H}_*^i(\E^\bullet)$ the $i$-th hypercohomology of
$\E^\bullet$.

\begin{prop}\label{P:spseq}
In the above situation, suppose that $\dim(Y) \le 0$. Then
\[
\mathbf{H}_*^0(\E^\bullet)=H^0_*(\mathcal{O}_Y)
\]
\[
\mathrm{Im}(\mathbf{d}_{r+1}) = H^0_{\mathbf{m}}(R), \quad
A_\infty^{00}=R/H^0_{\mathbf{m}}(R), \quad
A_\infty^{-rr}=H^1_{\mathbf{m}}(R)
\]
 and the exact sequence of edge homomorphisms
\[
\xymatrix{ 0\ar[r]&
A_\infty^{00}\ar[r]&\mathbf{H}_*^0(\E^\bullet)\ar[r]&
A_\infty^{-rr}\ar[r]& 0} \]
 coincides with the sequence:
\begin{equation}\label{E:edge}
\xymatrix{ 0\ar[r]&R/H^0_{\mathbf{m}}(R)\ar[r]&H^0_*(\mathcal{O}_Y)\ar[r]
&H^1_{\mathbf{m}}(R)\ar[r]&0}
\end{equation}
\end{prop}

\proof Let $x \in \P^r$. Then
\[
 \mathrm{depth}_x(\I_Y)\begin{cases}\ge r& \text{if $x \in Y$} \\ $=r+1$&
 \text{otherwise} \end{cases}
 \]
Therefore, by \cite{dE95}, Thm. 17.4 p. 424,  $(\H^q)_x= 0$ if
$q\le -2$ for all $x \in \P^r$, and $(\H^{-1})_x=0$ if $x \notin
Y$. Therefore $\H^q=0$ if $q \le -2$ and $\H^{-1}$ is supported on
$Y$. It follows that $H^p(\H^{-1})=0$
for all $p>0$.
 Now we decompose $\E^\bullet$ into short  exact sequences
of sheaves as follows:
\begin{equation}\label{E:E}
\xymatrix{ 0\ar[r]&\E_{-r-1}\ar[r]& \E_{-r}\ar[r]&I_{-r+1}\ar[r]&0
}\end{equation}
%\begin{equation}\label{E:H}
%\xymatrix{ &0\ar[r] & I_{-r+1}\ar[r]& K_{-r+1}\ar[r]&
%\H^{-r+1}\ar[r] &0} \end{equation}
\begin{equation}\label{E:K}
\xymatrix{ 0\ar[r]&
I_{-r+1}\ar[r]&\E_{-r+1}\ar[r]&I_{-r+2}\ar[r]&0 }
\end{equation}
etc., up to:
\[
\xymatrix{ 0 \ar[r]&I_{-2}\ar[r]&\E_{-2}\ar[r]&I_{-1}\ar[r]&0}
\]
\begin{equation}\label{E:K1}
\xymatrix{ &0\ar[r] & I_{-1}\ar[r]& K_{-1}\ar[r]& \H^{-1}\ar[r]
&0}
\end{equation}
\begin{equation}\label{E:Y}
\xymatrix{ 0\ar[r]& K_{-1}\ar[r]&\E_{-1}\ar[r]&\I_Y\ar[r]&0 }
\end{equation}
  The map $\mathbf{d}_{r+1}$ is obtained from a
diagram chasing out of these sequences. Since the $\E_i$'s are
direct sums of $\mathcal{O}(k)$'s, from (\ref{E:E}) and comparing with
(\ref{Err}) we deduce
\[
A_{r+1}^{-r-1r} \cong H^{r-1}_*(I_{-r+1})
\]
and from (\ref{E:K}), etc,   we have  isomorphisms
\[
A_{r+1}^{-r-1r} \cong H^{r-1}_*(I_{-r+1})\cong \cdots \cong
H^1_*(I_{-1})
\]

 Now we use (\ref{E:K1}) and we obtain   a surjective map:
\[
\xymatrix{ H^1_*(I_{-1}) \ar[r]& H^1_*(K_{-1})\ar[r]&0}
\]
But from sequence (\ref{E:Y}) it follows that
\[
H^1_*(K_{-1})= H^0_{\mathbf{m}}(R) 
\]
 and this proves that
Im$(\mathbf{d}_{r+1})=H^0_{\mathbf{m}}(R)$. Therefore it also follows that
\[
A_\infty^{00}= A^{00}_{r+1}/\mathrm{Im}(\mathbf{d}_{r+1})=
R/H^0_{\mathbf{m}}(R)
\]
Now observe that $A^{-rr}_\infty = H^r_*(I_{-r+1})$. A diagram
chasing similar to the previous one shows that
\[
H^r_*(I_{-r+1}) \cong H^1_*(\I_Y)
\]
Since $H^1_*(\I_Y)= H^1_{\mathbf{m}}(R)$ we obtain the
identification $A^{-rr}_\infty =H^1_{\mathbf{m}}(R)$.

Noting that  the $B$-spectral sequence degenerates at $B_2$,   we
get in particular that
\[
\mathbf{H}_*^0(\E^\bullet)= H^0_*(\H^0(\E^\bullet))=H^0_*(\mathcal{O}_Y)
\]
Therefore the edge exact sequence is (\ref{E:edge}). \qed

As a consequence we can now derive the following:

\begin{theorem}\label{C:spseq}
Let $I=(f_0,\dots,f_r)$ with deg$(f_j)=d_j$, $R=P/I$ and $Y = \mathrm{Proj}(R)$. 
Assume that $\mathrm{dim}(Y)\le 0$. Then
 there is a natural isomorphism:
 \[
 H^0_{\mathbf{m}}(R)\cong [H^0_{\mathbf{m}}(R)(\sigma)]^\vee
 \]
 where $\sigma = \sum_{j=0}^r d_j -r-1$. Therefore we have
 natural isomorphisms
 \[
 H^0_{\mathbf{m}}(R)_k \cong H^0_{\mathbf{m}}(R)_{\sigma-k}^\vee, \quad 0 \le k \le \sigma
 \]

\end{theorem}

\proof  
The surjective map: \[
\mathbf{d}_{r+1}:\xymatrix{[R(\sigma)]^\vee\ar[r]&H^0_{\mathbf{m}}(R)}
\]
  dualizes as an injective map:
\[
\mathbf{d}^\vee_{r+1}:\xymatrix{ H^0_{\mathbf{m}}(R)^\vee\ar[r]&R(\sigma)}
\]
whose image must be contained in $H^0_{\mathbf{m}}(R)(\sigma)$
because it consists of elements which are killed by
$\mathbf{m}^{\sigma+1}$. But then
$\mathrm{Im}(\mathbf{d}^\vee_{r+1})=H^0_{\mathbf{m}}(R)(\sigma)$   because
$H^0_{\mathbf{m}}(R)^\vee$ and $H^0_{\mathbf{m}}(R)(\sigma)$ have
the same dimension as vector spaces. \qed

\begin{remark}\rm\label{R:spseq1}
As already stated in the Introduction, Corollary \ref{C:spseq} is a special case of \cite{vSW11}, Theorem 4.7.
The case $Y = \emptyset$ of course corresponds
to the situation when the elements $f_0,\dots, f_r$ form a regular sequence, and this happens if and only
if $\H^q(\E^\bullet)=0$ for all $q$. In this case the
hypercohomology $\mathbf{H}_*^\bullet(\E^\bullet)$ is zero in all
dimensions, because the $B_2$-spectral sequence is zero. It
follows that the map:
\[
\mathbf{d}_{r+1}:A_{r+1}^{-r-1r}\longrightarrow A_{r+1}^{00}
\]
is an isomorphism, which means that we have an isomorphism
$[R(\sigma)]^\vee \cong R$. This is  the well known duality theorem of
Macaulay for Gorenstein artinian algebras (\cite{cV03}, Th. II6.19, p. 172).
\end{remark}

As a special case of Theorem \ref{C:spseq} we obtain the following (see also \cite{DS12}, Theorem 1):

\begin{theorem}\label{T:selfdu}
Assume that the hypersurface $X$ has at most isolated
    singularities.  Then:
    \[
H^0_{\mathbf{m}}(R(f))_k \cong
H^0_{\mathbf{m}}(R(f))_{\sigma-k}^\vee, \quad 0 \le k \le \sigma
 \]
where $\sigma = (r+1)(d-2)$.
\end{theorem}

\begin{remark}\rm\label{R:spseq2}
In case $X$ is nonsingular the jacobian ring $R=R(f)$ is Gorenstein artinian with socle in 
degree $\sigma$. The self duality of $R(f)$ is induced by a pairing
\[
R_k \times R_{\sigma-k} \longrightarrow R_\sigma \cong \mathbb{C}
\]
where the first map is induced by multiplication of polynomials and the last isomorphism
is obtained from the trace map for local duality.

 \end{remark}

\begin{corollary}\label{C:log1}
Assume that $X$ has only isolated singularities. Then there are natural isomorphisms:
\[
H^1(T\langle X \rangle(-d+k)) \cong H^1(T\langle X \rangle(\sigma-d-k))^\vee
\]
for all $k$. 
\end{corollary}

\proof Use (\ref{E:log5}) and Theorem (\ref{T:selfdu}). \qed 

Observe that, in case  the hypersurface $X$ is singular with isolated singularities and $r \ge 3$, the sheaf 
$T\langle X \rangle(-d)$ is reflexive of rank $r$ but not locally free (see \cite{kS80}). Therefore the duality statement of Corollary
\ref{C:log1} is not a consequence of standard properties of locally free sheaves.  

On the other hand if $r=2$
then $T\langle X \rangle(-d)$ is locally free of rank two and its first Chern class is given by:
\[
c_1(T\langle X \rangle(-d)) = 3-3d
\]
Then  Corollary
\ref{C:log1} follows directly from the  straightforward fact that for every locally free sheaf $E$ of rank two on $\P^2$ 
 we have 
\[
H^1(E(k)) \cong H^1(E(\sigma-k))^\vee
\]
where $\sigma=-c_1(E)-3$. 

It is not clear how far one can go relaxing the hypothesis of Theorem \ref{T:selfdu}, as the next two examples show.

 \begin{example}\label{Ex:ruledsur} \rm
The \emph{ruled cubic surface} $X \subset \P^3$ has equation
\[
XT^2-YZ^2=0
\]
and is singular along the line $T=Z=0$. The local cohomology has
only one non-zero term in degree 2, and:
\[
h^0_{\mathbf{m}}(R(f))_2= 1
\]
Since $\sigma=4$, the symmetry condition $H^0_{\mathbf{m}}(R(f))_k
\cong H^0_{\mathbf{m}}(R(f))_{\sigma-k}$ is fullfilled even though $X$ doesn't   satisfy the hypothesis of
Theorem \ref{T:selfdu}.
\end{example}

\begin{example}\label{Ex:quarsur} \rm
A \emph{quartic surface with a double conic}  $X \subset \P^3$ has equation:
\[
 (ZT-XY)^2 + (X+Y+Z+T)^2(X^2+Y^2+Z^2+T^2)=0
\]
The table of its local cohomology dimensions is:

\begin{tabular}{|c|c|}
  \hline
   $j$ & $h^0_{\mathbf{m}}(R)_j$  \\
  \hline
  0 & 0   \\
  1 & 0  \\
  2 & 1  \\
  3 & 4   \\
  4 & 5   \\
  5 & 1    \\
  6 & 0   \\
  7 & 0   \\
  8 & 0   \\
  \hline
\end{tabular}

Since $\sigma = 10$, we see that   self-duality does not hold in this case.
\end{example}

%%%%%%%%%%%%%%%%%%%%%%%%%%%%%%%%%%%%%%%%%%%%%%%%%%%%%

 %%%%%%%%%%%%%%%%%%%%%%%%%%%%%%%%%%%%%%%%%%%%%%%%%%%%%

\section{Logarithmic differentials}\label{S:logar}

	Let's restrict for a moment to the case when our $X \subset \P^r$ of degree $d$ is nonsingular.
	Then Griffiths' Theorem identifies:
	\begin{equation}\label{E:griff1}
	\bigoplus_{p=1}^r H^{r-p,p-1}(X)_0  = \bigoplus_{p=1}^r H^1(T\langle X \rangle (K_{\P^r}+(p-1)X)  
	\end{equation}
	thanks to  Proposition \ref{P:log1}, which identifies 
	\[
	\bigoplus_{p=1}^r H^1(T\langle X \rangle (K_{\P^r}+(p-1)X) = \bigoplus_{p=1}^r R(f)_{pd-r-1}
		\]
		 The right hand side of (\ref{E:griff1}) is well defined if 
		$X$ is just a reduced hypersurface in a projective manifold $Z$ of dimension $r$, after replacing $\P^r$ with $Z$. 
		In such a situation  it is convenient to consider, together with $T_Z\langle
X\rangle$, the \emph{sheaves of logarithmic differentials along
$X$} which are defined as follows:
  \[
  \Omega_Z^k(\mathrm{log}X):= \{\omega\in \Omega_Z^k(X):
  d\omega\in\Omega_Z^{k+1}(X)\}, \quad k=0,\dots, r
\]
In particular $\Omega_Z^0(\mathrm{log}X)= \mathcal{O}_Z$ and
$\Omega_Z^r(\mathrm{log}X)= K_Z+X$. For $k\ne 0,r$ these
sheaves are not locally free in general.  For $k=1$ one has:
\[
\Omega_Z^1(\mathrm{log}X) := Hom_Z(T_Z\langle
X\rangle,\mathcal{O}_Z)
\]
and this sheaf is reflexive (\cite{kS80}, n. 1.7). By definition we have inclusions
\[
\Omega^k_Z \subset \Omega_Z^k(\mathrm{log}X) \subset \Omega^k_Z(X)
\]
which in turn induce  the  inclusions:
\begin{equation}\label{E:logincl}
  \Omega^k_Z(\mathrm{log}X)(-X) \subset \Omega^k_Z\subset
\Omega^k_Z(\mathrm{log}X)  
\end{equation}
We collect in the following Lemmas the properties we need about
the sheaves of logarithmic differentials.

\begin{lemma}\label{L:hodge1}
The following conditions are equivalent:

\begin{itemize}

\item[(i)] $T_Z\langle X\rangle$ is locally free.

\item[(ii)] \[ \Omega_Z^k(\mathrm{log}X)= \bigwedge^k
\Omega_Z^1(\mathrm{log}X)
\]
for all $k=1,\dots, r$.

    \item[(iii)] $\bigwedge^r
\Omega_Z^1(\mathrm{log}X) = \Omega_Z^r(\mathrm{log}X)
(=K_Z+X)$
\end{itemize}
If the above conditions are satisfied then we have a canonical
identification:
    \begin{equation}\label{E:dualog}
 T_Z\langle
X\rangle(K_Z+X) = \Omega_Z^{r-1}(\mathrm{log}X)
\end{equation}
\end{lemma}

\proof  The equivalence of the conditions stated is Theorem 1.8 of
\cite{kS80}.  From (iii) we obtain $c_1(T_Z\langle
X\rangle)=-(K_Z+X)$. Therefore:
\begin{align}
T_Z\langle X\rangle(K_Z+X) &=T_Z\langle X\rangle
c_1(T_Z\langle X\rangle^\vee) \notag \\
&=\bigwedge^{r-1}T_Z\langle X\rangle^\vee \notag \\
& =\bigwedge^{r-1}\Omega_Z^1(\mathrm{log}X)\notag
\\
\text{by (i)}&=\Omega_Z^{r-1}(\mathrm{log}X)\notag
\end{align} \qed

The following are examples   such that $T_Z\langle X \rangle$ is locally free
(see \cite{kS80}):
\begin{itemize}
    \item  $X$ nonsingular.
    \item  $Z$ is a surface ($r=2$).
    \item $X$ has normal crossing singularities at every point (it is a normal crossing divisor). Recall
    that this means that for each $x \in X$ the local ring $\mathcal{O}_{X,x}$ is
    formally, or etale, equivalent to $\mathcal{O}_{Z,x}/(t_1\cdots t_k)$ for some $1 \le k \le r-1$, where
    $t_1, \dots, t_k$ are part of a local system of coordinates.
     \end{itemize}

     Recall that $X\subset Z$ is a \emph{strictly normal crossing divisor}
     if it is a normal crossing divisor whose irreducible
     components $X_1, \dots, X_s$ are nonsingular.

\begin{lemma}\label{L:hodge2}
Assume that $X = X_1 \cup \dots \cup X_s \subset Z$ is a strictly normal crossing divisor. Denote by
$\widehat{X}_1 = X_2\cap \dots \cap X_s$, and by $Y_1 = X_1 \cap \widehat{X}_1$.  Then there
are exact sequences, for $a=1, \dots, r=\dim (Z)$:
\begin{align}
    \xymatrix{
    0 \ar[r] & \Omega^1_Z \ar[r] & \Omega_Z^1(\mathrm{log} X)\ar[r] & \bigoplus_{i=1}^s \mathcal{O}_{X_i} \ar[r] & 0}\label{E:hv1} \\
	\xymatrix{
    0 \ar[r] & \Omega^a_Z(\mathrm{log} \widehat{X}_1) \ar[r] & \Omega^a_Z(\mathrm{log} X) \ar[r]^-R & \Omega^{a-1}_{X_1}(\mathrm{log} Y_1)\ar[r]& 0} \label{E:hv2} \\
		\xymatrix{
		 0 \ar[r] &  \Omega^a_Z(\mathrm{log} X)(-X_1) \ar[r] & \Omega^a_Z(\mathrm{log} \widehat{X}_1) \ar[r] & \Omega^a_{X_1}(\mathrm{log} Y_1)\ar[r]& 0} \label{E:hv3}
\end{align}
where $R$ is the \emph{residue operator}.
 \end{lemma}

\proof   see  \cite{EV92}, \S 2.3. \qed

Note that, by twisting (\ref{E:hv3}) by $\mathcal{O}_Z(-\widehat{X}_1)$ we obtain the following exact sequence:
\begin{equation}\label{E:hv4}
\xymatrix{
		 0 \ar[r] &  \Omega^a_Z(\mathrm{log} X)(-X) \ar[r] & \Omega^a_Z(\mathrm{log} \widehat{X}_1)(-\widehat{X}_1) \ar[r] & \Omega^a_{X_1}(\mathrm{log} Y_1)(-Y_1)\ar[r]& 0}
 \end{equation}
For future reference it is worth emphasizing that when $X=X_1$ is irreducible and nonsingular then the   sequences (\ref{E:hv2}) and (\ref{E:hv4})   become respectively:
 \begin{align}
    \xymatrix{
    0 \ar[r] &
    \Omega^a_Z\ar[r]&\Omega^a_Z(\mathrm{log}X)\ar[r]^-R& \Omega^{a-1}_X\ar[r]&0}\label{E:logseq1}\\
\xymatrix{
    0 \ar[r] &\Omega^a_Z(\mathrm{log}X)(-X)\ar[r] &
    \Omega^a_Z\ar[r]&
    \Omega^a_X\ar[r]&0}\label{E:logseq2}
\end{align}

		\begin{lemma}\label{L:hodge3}
		Assume that $X\subset Z$ is an irreducible and nonsingular divisor.
		For each $k=0, \dots, r-1$ consider the composition:
		\[
		\lambda:\xymatrix{
H^k(\Omega^k_X)\ar[r]^-\delta&H^{k+1}(\Omega^{k+1}_Z)\ar[r]^-{\nu_{k+1}^*}&
H^{k+1}(\Omega^{k+1}_X)}
\]
where $\delta$ is a coboundary map of the sequence (\ref{E:logseq1}) and
$\nu_{k+1}^*$ is induced by the second homomorphism in
the sequence (\ref{E:logseq2}).  Then $\lambda$ is the map
defined by the Lefschetz operator corresponding to the Kahler
metric on $X$ associated to   
$\mathcal{O}_X(X)$.
\end{lemma}
		
\proof The Lefschetz operator $L:H^k(X,\mathbb{C})\longrightarrow H^{k+2}(X,\mathbb{C})$
is the composition:
\[
 	\xymatrix{
		H^k(X,\mathbb{C})\ar[r]^-{\gamma}&H^{k+2}(Z,\mathbb{C}) \ar[r]^-{\nu^*}&H^{k+2}(X,\mathbb{C})
		}
		\]
	where $\gamma$ is the Gysin map and $\nu^*$ is induced by the inclusion 
	\[
	\xymatrix{
	 X \ar@{^(->}[r]^-{\nu}& Z
	}\]
	(\cite{cV03}, v. II, (2.11) p. 57).	
	Moreover $\gamma$ is the cokernel of the  map:
	\[
	\rho:H^{k+1}(U,\mathbb{C}) \longrightarrow H^k(X,\mathbb{C})
	\]
	induced by the residue operator, where $U=Z\setminus X$.  
	More precisely,   we have an  isomorphism (\cite{cV03},   Corollary I.8.19 p. 198)
	\[
	H^\bullet(U,\mathbb{C}) \cong \mathbb{H}^\bullet(\Omega^\bullet(\mathrm{log} X))
	\]
	(where $\mathbb{H}$ denotes hypercohomology) and the   map $\rho$ is induced by the  residue operators $R$ of the exact sequences 
	(\ref{E:logseq1}). Therefore  the restriction of $\gamma$ to $H^k(\Omega^k_{X})$ is identified with $\delta$ (see \cite{cV03}, Prop. I.8.34 p. 210).
	On the other hand   $\nu_{k+1}^*$ is the restriction of $\nu^*$
	to $H^{k+1}(\Omega^{k+1}_Z)$.  \qed

	%%%%%%%%%%%%%%%%%%%%%%%%%%%%%%%%%%%%%%%%%%%%%%%%%%%%%%%%%%%%%%%%%%%%%
	
	\section{Local cohomology and Hodge theory}\label{S:hodge}

	We now come back to the original situation of a reduced hypersurface $X =V(f)\subset \P^r$ of degree $d$. By Proposition \ref{P:log1}
	for $p=1, \dots, r$ we can identify 
	\begin{equation}\label{E:ncd2}
	H^0_{\mathbf{m}}(R(f))_{pd-r-1}= H^1_*(T\langle X \rangle(K_{\P^r}+(p-1)X))
		\end{equation}
		Moreover, \emph{if $T\langle X \rangle$ is locally free} then, by Lemma \ref{L:hodge1}, we also have:
	\begin{equation}\label{E:ncd3}
	H^0_{\mathbf{m}}(R(f))_{pd-r-1} =  H^1(\Omega^{r-1}(\mathrm{log}X)(p-2)X)
	\end{equation}
	Our first result is the following:

	\begin{theorem}\label{T:ncd2}
	Assume that   $X\subset \P^r$ is a strictly normal crossing hypersurface, with irreducible components $X_1, \dots, X_s$. 
	 Then we have:
	\[
	H^0_{\mathbf{m}}(R(f))_{d-r-1} \cong \bigoplus_{i=1}^s H^0(X_i,\Omega^{r-1}_{X_i})
	\]

	\end{theorem}
	
	\proof Since $T\langle X\rangle$ is locally free we have the identification (\ref{E:ncd3}) for $p=1$:
	\[
	H^0_{\mathbf{m}}(R(f))_{d-r-1} = H^1(\Omega^{r-1}(\mathrm{log}X)(-X))
	\]
	Assume first $r \ge 3$. Consider the exact sequence (\ref{E:hv4}) for $a=r-1$.  Since 
	\[
	h^0(\P^r,\Omega^{r-1}_{\P^r}(\mathrm{log} \widehat{X}_1)(-\widehat{X}_1)) =0
	\]
	we obtain the exact sequence:
	\[
	\xymatrix{
	0 \ar[r] & H^0(\Omega_{X_1}^{r-1})\ar[r]& H^0_{\mathbf{m}}(R(f))_{d-r-1} \ar[r]&
	H^1(\Omega^{r-1}_{\P^r}\mathrm{log} \widehat{X}_1)(-\widehat{X}_1))\ar[r]&0
	}
	\]
	where the zero on the right is $H^1(\Omega_{X_1}^{r-1})$. Now the conclusion follows by induction on $s$. 
	
	If $r=2$  and $s=1$ use (\ref{E:logseq2}) and Lemma \ref{L:hodge3}. If $s \ge 2$ use (\ref{E:hv4}) and induction.
	\qed

	\begin{corollary}\label{C:ncd1}
	Let $X=X_1+ \cdots+ X_s\subset \P^2$ be a strictly normal crossing plane curve.   Then
	\[
	H^0_{\mathbf{m}}(R(f))_{d-3} \cong \bigoplus_{i=1}^s H^0(X_i,\omega_{X_i}), \quad H^0_{\mathbf{m}}(R(f))_{2d-3} \cong \bigoplus_{i=1}^s H^1(X_i,\mathcal{O}_{X_i})
	\]
	\end{corollary}
	
	\proof It follows from the theorem, from the self duality theorem \ref{T:selfdu}, and Serre duality applied to each component $X_i$.  \qed

	When $r \ge 3$ the relation  between  the other graded pieces $H^0_{\mathbf{m}}(R(f))_{pd-r-1}$, $p=2,\dots, r$, of the local cohomology 
	and the primitive middle cohomology of the components of $X$ is more complicated because the  intersections of the components contribute non-trivially. 
	As an example we compute the   dimension of the middle 
	term in the case $r=3$.
	
	\begin{theorem}\label{T:ncd3}
	Let $X=X_1+ \cdots+ X_s \subset \P^3$ be a strictly normal crossing surface, whose components have degrees $d_1, \dots, d_s$ respectively. Then:
	\[
	h^0_{\mathbf{m}}(R(f))_{2d-4} = \sum_{i=1}^s \dim[H^{1,1}(X_i)_0] + \sum_{1\le i< j \le s} g(X_i\cap X_j)
	\]
	where $g(X_i\cap X_j) = \frac{1}{2}d_id_j(d_i+d_j-4)+1$ is the genus of the curve $X_i\cap X_j$.
	\end{theorem}
	
	\proof By induction on $s$.  If $s=1$ the formula is true by Griffiths' Theorem. Assume $s \ge 2$. Then 
	$H^0_{\mathbf{m}}(R(f))_{2d-4} = H^1(\Omega^2_{\P^3}(\mathrm{log} X))$, by (\ref{E:ncd3}).
	We let 
	\begin{align}
	 \widehat{X}_1&=X_2+\cdots + X_s \notag \\
	Y_1 &= X_1 \cap (X_2+\cdots + X_s) \notag \\
	  \widehat{Y}_1&= X_1 \cap (X_3+\cdots + X_s)\notag
		\end{align}
		 We have the following diagram of exact sequences:
	\begin{equation}\label{E:ncd6}
	\xymatrix{
	&&&0\ar[d] \\ &&&\Omega^1_{X_1}(\mathrm{log} \widehat{Y}_1) \ar[d] \\
	0 \ar[r] & \Omega^2_{\P^3}(\mathrm{log} \widehat{X}_1) \ar[r] & \Omega^2_{\P^3}(\mathrm{log} X) \ar[r] & \Omega^1_{X_1}(\mathrm{log} Y_1)\ar[r]\ar[d]& 0 \\
	&&& \mathcal{O}_{X_1\cap X_2}\ar[d] \\
	&&&0} 
	\end{equation}
	We claim the following:
	\begin{itemize}
		\item[(a)] $h^0(\Omega^1_{X_1}(\mathrm{log} Y_1))=s-2$.
		\item[(b)] $H^2(\Omega^2_{\P^3}(\mathrm{log} \widehat{X}_1))=0$.
		\item[(c)] $H^2(\Omega^1_{X_1}(\mathrm{log} \widehat{Y}_1)) =0$.
	\end{itemize}
	Assume that (a),(b),(c) are proved. Then from the above diagram we deduce  the exact sequence:
	\begin{equation}\label{E:ncd4}
	\xymatrix{
	0 \ar[r] & H^1(\Omega^2_{\P^3}(\mathrm{log} \widehat{X}_1))  \ar[r] & H^1(\Omega^2_{\P^3}(\mathrm{log} X)) \ar[r] & H^1(\Omega^1_{X_1}(\mathrm{log} Y_1)) \ar[r] &0}
	\end{equation}
	The   term on the right in (\ref{E:ncd4}) can be computed using the vertical exact sequence of   diagram (\ref{E:ncd6}).
	 Assume first that $s=2$. In this case $Y_1=X_1\cap X_2$ and recalling (a) we obtain:
	\begin{small}\[
	\xymatrix{
	0 \ar[r]&  H^0(\mathcal{O}_{X_1\cap X_2})\ar[r]&H^1(\Omega^1_{X_1})\ar[r] & H^1(\Omega^1_{X_1}(\mathrm{log} (X_1\cap X_2))\ar[r]&H^1(\mathcal{O}_{X_1\cap X_2})\ar[r]&0 }
	\]\end{small}
	whence:
	\[
	H^1(\Omega^1_{X_1}(\mathrm{log} (X_1\cap X_2)) = \dim[H^{1,1}(X_1)_0]+ g(X_1\cap X_2)
	\]
	If $s\ge 3$ then the map 
	\[
	H^0(\mathcal{O}_{X_1\cap X_2})\longrightarrow H^1(\Omega^1_{X_1}(\mathrm{log} \widehat{Y}_1))
	\]
	is zero by (a). Therefore, applying induction, from the vertical exact sequence of   diagram (\ref{E:ncd6}) we deduce:
	\begin{align}                                                                       
	\dim[H^1(\Omega^1_{X_1}(\mathrm{log} Y_1))]& = \dim[H^1(\Omega^1_{X_1}(\mathrm{log} \widehat{Y}_1))]+ g(X_1\cap X_2) \notag \\
	&= \dim[H^{1,1}(X_1)_0]+ \sum_{i=3}^s g(X_1\cap X_i)+g(X_1\cap X_2) \notag \\
	&= \dim[H^{1,1}(X_1)_0]+ \sum_{i=2}^s g(X_1\cap X_i) \notag
	\end{align}
	By induction we have:
	\[
	h^1(\Omega^2_{\P^3}(\mathrm{log} \widehat{X}_1)) = \sum_{i=2}^s \dim[H^{1,1}(X_i)_0] + \sum_{2\le i< j \le s} g(X_i\cap X_j)
	\]
	Therefore, putting all these computations together the claimed expression for $h^0_{\mathbf{m}}(R(f))_{2d-4}$ follows. We still have to prove (a),(b) and (c).
	
	\separation
	
	\emph{Proof of (a).}  Use the exact sequence   (\ref{E:hv1}) with $Z=X_1$ and $X=Y_1$, and the fact that the image of the coboundary map
	is the space generated by the Chern classes of $X_1\cap X_2, \dots, X_1 \cap X_s$, which is 1-dimensional.
	
	\separation
	
	\emph{Proof of (c).}  Use the vertical sequence in (\ref{E:ncd6}) and induction on $s \ge 2$.
	
	\separation
	
	\emph{Proof of (b).} Assume $s=1$. The    map $H^1(\Omega_{X_1}^1)\longrightarrow H^2(\Omega_{\P^3}^2)$ coming from the sequence 
	\[
	\xymatrix{
	0 \ar[r] & \Omega^2_{\P^3}\ar[r] & \Omega^2_{\P^3}(\mathrm{log} X_1)\ar[r] & \Omega_{X_1}^1 \ar[r] & 0}
	\]
	is surjective (this follows from Lemma \ref{L:hodge3} and Hodge theory). Therefore since $H^2(\Omega^1_{X_1})=0$, it follows that
	$H^2(\Omega^2_{\P^3}(\mathrm{log} X_1))=0$. The general case of (b) now follows by induction, from (c)  and from the exact row in (\ref{E:ncd6}).
	\qed
	
	We give a few examples illustrating these results.

\begin{example}\label{Ex:redqua1} \rm
Let $f=X_0(X_0^3+X_1^3+X_2^3)$.  Then $X=V(f)=\Lambda\cup C\subset
\P^2$ is a \emph{strictly normal crossing reducible plane quartic},
consisting of a line $\Lambda$ and a nonsingular cubic $C$. One
computes that:
\begin{align}
 H^0_{\mathbf{m}}(R)_1&\cong \mathbb{C}\cong
H^{1,0}(X')=H^{1,0}(C)\notag \\
H^0_{\mathbf{m}}(R)_5&\cong\mathbb{C}\cong
H^{0,1}(X')=H^{0,1}(C)\notag
\end{align}
The complete table is:

\medskip

\begin{tabular}{|c|c|c|c|}
  \hline
  % after \\: \hline or \cline{col1-col2} \cline{col3-col4} ...
  $j$ & $h^0_{\mathbf{m}}(R)_j$ & $\dim(R(f)_j)$ & $h^0(\mathcal{O}_{\mathrm{Sing}(X)}(j))$ \\
  \hline
  0 & 0 & 1 & 1 \\
  1 & 1 & 3 & 2 \\
  2 & 3 & 6 & 3 \\
  3 & 4 & 7 & 3 \\
  4 & 3 & 6 & 3 \\
  5 & 1 & 4 & 3 \\
  6 & 0 & 3 & 3 \\
  7 & 0 & 3 & 3 \\
  8 & 0 & 3 & 3 \\
  \hline
\end{tabular}
\end{example}

The conclusion of Corollary \ref{C:ncd1} fails even in the simplest cases if one weakens the assumptions about the singularities of $X$, 
as the following two examples show.

\begin{example}\label{Ex:cusqua} \rm
A \emph{1-cuspidal plane quartic} $f =
X_0^2X_1^2+X_1^2X_2^2+X_1^4+X_2^4$.  Here the table is:

\medskip

\begin{tabular}{|c|c|c|c|}
  \hline
  % after \\: \hline or \cline{col1-col2} \cline{col3-col4} ...
  $j$ & $h^0_{\mathbf{m}}(R)_j$ & $\dim(R(f)_j)$ & $h^0(\mathcal{O}_{\mathrm{Sing}(X)}(j))$ \\
  \hline
  0 & 0 & 1 & 1 \\
  1 & 1 & 3 & 2 \\
  2 & 4 & 6 & 2 \\
  3 & 5 & 7 & 2 \\
  4 & 4 & 6 & 2 \\
  5 & 1 & 3 & 2 \\
  6 & 0 & 1 & 2 \\
  7 & 0 & 1 & 2 \\
  8 & 0 & 1 & 2 \\
  \hline
\end{tabular}

\medskip

 Then  $X'$ has genus two, has self-dual local cohomology but 
 $h^0_{\mathbf{m}}(R)_1=1=h^0_{\mathbf{m}}(R)_5 < 2$.  
\end{example}

\begin{example}\label{Ex:redqua2} \rm
A \emph{reducible plane quartic} consisting of a nonsingular cubic and
of an inflectional tangent:
\[
f = X_0(X_0^2X_1+X_0X_1^2+X_2^3)
\]
In this case the table is:

\medskip

\begin{tabular}{|c|c|c|c|}
  \hline
  % after \\: \hline or \cline{col1-col2} \cline{col3-col4} ...
  $j$ & $h^0_{\mathbf{m}}(R)_j$ & $\dim(R(f)_j)$ & $h^0(\mathcal{O}_{\mathrm{Sing}(X)}(j))$ \\
  \hline
  0 & 0 & 1 & 1 \\
  1 & 0 & 3 & 3 \\
  2 & 1 & 6 & 5 \\
  3 & 2 & 7 & 5 \\
  4 & 1 & 6 & 5 \\
  5 & 0 & 5 & 5 \\
  6 & 0 & 5 & 5 \\
  7 & 0 & 5 & 5 \\
  8 & 0 & 5 & 5 \\
  \hline
\end{tabular}
\end{example}

\begin{example}\label{Ex:redsurf} \rm
 \emph{A strictly normal crossing quintic surface.}
(This example has been kindly suggested by A. Dimca).
As an illustration of Theorem \ref{T:ncd3} consider $X=V(f) \subset \P^3$, where 
\[
f(X_0,\dots,X_3) = (X_0^2+X_1^2+X_2^2+X_3^2)(X_0^3+X_1^3+X_2^3+X_3^3)
\]
Then $X = X_1+X_2$ is the union of a quadric and a cubic, and $C = X_1 \cap X_2$ is a canonical curve (of genus $4$).
The table of local cohomology is:
\medskip

\begin{tabular}{|c|c|c|c|}
  \hline
  % after \\: \hline or \cline{col1-col2} \cline{col3-col4} ...
  $j$ & $h^0_{\mathbf{m}}(R)_j$ & $\dim(R(f)_j)$ & $h^0(\mathcal{O}_C(j))$ \\
  \hline
  0 & 0 & 1 & 1 \\
  1 & 0 & 4 & 4 \\
  2 & 1 & 10 & 9 \\
  3 & 5 & 20 & 15 \\
  4 & 10 & 31 & 21 \\
  5 & 13 & 40 & 27 \\
  6 & 11 & 44 & 33 \\
  7 & 5 & 44 & 39 \\
  8 & 1 & 46 & 45 \\
	9 &0&51&51\\
	10 &0&57&57\\
	11 &0&63&63\\
	12 &0&69&69 \\
	13 &0&75&75\\
	
  \hline
\end{tabular}

\medskip

Note that 
\[
h^0_{\mathbf{m}}(R)_6 = 11 = (2-1)+(7-1)+4 = h^{1,1}(X_1)+h^{1,1}(X_2)+g(C)
\]
 as expected.
\end{example}

%%%%%%%%%%%%%%%%%%%%%%%%%%%%%%%%%%%%%%%%%%%%%%%%%%%%%%%

\section{Torelli-type questions}\label{S:torelli}

Following a terminology introduced in \cite{DK93}, a reduced hypersurface $X\subset \P^r$ is called \emph{Torelli}
 in the sense of Dolgachev-Kapranov if it can be reconstructed from the sheaf $T\langle X\rangle$. In their paper 
\cite{DK93}  they studied the Torelli property of normal crossing arrangements of hyperplanes. Their main result has been 
later improved by Vall\`es in \cite{jV00}.  For arbitrary arrangements of hyperplanes the Torelli problem has been settled in
\cite{FMV}.
In \cite{UY09} it is proved that a smooth hypersurface is Torelli if and only if it is not of
Sebastiani-Thom type.
 E. Angelini \cite{eA13} studied certain normal crossing configurations of smooth hypersurfaces proving that they are Torelli in several cases.

We want to consider a different reconstruction problem, namely we ask:

\begin{itemize}
	\item[] \textbf{Question:} Under which circumstances can $X$ be reconstructed from 
	\[
	H^0_{\mathbf{m}}(R(f)) \cong H^1_*(T\langle X \rangle(-d))
	\]
\end{itemize}

In the nonsingular case this   is merely the question of reconstructability of $X$ from its jacobian ring.
  This question has been considered extensively in the literature, even in the singular case.
The typical result one would like to generalize is the following:

\begin{theorem}\label{T:donagi}
\begin{itemize}
	 \item[(i)]\cite{rD83} Let $f$ and $f'$ be homogeneous polynomials of degree $d$ defining  reduced hypersurfaces in $\P^r$. If $J(f)_d = J(f')_d$ 
	then $f$ and $f'$ are projectively equivalent.

	\item[(ii)] \cite{CG80} Let $f\in P$ be a generic polynomial of degree $d\ge 3$. 
	Then $f$ is determined by $J(f)_{d-1}$, up to a constant factor.
\end{itemize}
\end{theorem}

In this respect the following result is   relevant:

\begin{theorem}[\cite{gH64}]\label{T:hor2}
A locally free sheaf $\mathcal{F}$ of rank two on $\P^2$ can be reconstructed from the $P$-module $H^1_*(\mathcal{F})$.
\end{theorem}

Theorem \ref{T:hor2} suggests that, at least in $\P^2$, the reconstructability of $X$ from the module $H^1_*(T\langle X\rangle)$ 
is equivalent to the reconstructability of $X$ from the sheaf $T\langle X\rangle$. In fact we have the following: 

\begin{theorem}\label{T:tor1}
 A reduced plane curve is Torelli in the sense of Dolgachev-Kapranov if and only if it can be reconstructed from the local cohomology of its jacobian ring.  
 \end{theorem}

\proof It is an immediate consequence of Theorem \ref{T:hor2} and of the fact that $T\langle X\rangle$ is locally free for reduced plane curves. \qed

  Theorem \ref{T:tor1} of course applies to Torelli  arrangements of lines, that have been characterized as recalled above, and to  
	normal crossing arrangements of sufficiently many nonsingular curves of the same degree $n$ 
	(see \cite{eA13} for the precise statement). 
	Much less is known in the irreducible case, even for plane curves.  For partial results in this direction we refer the reader to \cite{DS14}.
	The Torelli property is related with freeness, that we are going to discuss next.

	%%%%%%%%%%%%%%%%%%%%%%%%%%%%%%%%%%%%%%%%%%%%%%%%%%%%%%

\section{Freeness}\label{S:free}

 According to Proposition \ref{P:log1} the vanishing of  $H^0_{\mathbf{m}}(R(f))$ is equivalent to that of 
$H^1_*(T\langle X\rangle)$ and it is a necessary condition for the freeness of $T\langle X\rangle$.
% since $T\langle X\rangle$ is reflexive, 
%Theorem \ref{T:hor1} implies that this occurs if and only if $T\langle X\rangle$ is  free. 
If $X$ is nonsingular then 
$H^0_{\mathbf{m}}(R(f))=R(f)$ is never zero, and therefore $T\langle X\rangle$ cannot be free. The same is true if
$\mathrm{Sing}(X)\ne \emptyset$ and has  codimension $\ge 2$ in $X$, because then  $T\langle X\rangle$ is not even locally free.

In   general     little seems to be known about the freeness of $T\langle X\rangle$, even in the case $r=2$. 
 We will mostly restrict to this case in the remaining of this section.

Look at the exact sequence: 
\begin{equation}\label{E:resol1}
\xymatrix{
0 \ar[r] & T\langle X\rangle (-1) \ar[r]& \mathcal{O}_{\P^2}^3 \ar[r]^-{\partial f} & \mathcal{J}_f(d-1) \ar[r] &0}
\end{equation}
Then
\[
c_1(T\langle X\rangle (-1)) = 1-d, \quad c_2(T\langle X\rangle (-1)) = (d-1)^2-t^1_X
\]
where $t^1_X=\dim(T^1_X)$. If    $T\langle X\rangle (-1)= \mathcal{O}(-a)\oplus \mathcal{O}(-b)$ is free then 
\begin{equation}\label{E:resol2}
a+b=d-1, \quad  ab=(d-1)^2-t^1_X
\end{equation}
They imply together that:
\begin{equation}\label{E:resol3}
a^2+ab+b^2 = t^1_X
\end{equation}
Observe also that, since under the restriction $a+b=d-1$ the product $ab$ attains its maximum when $(a,b)$ is balanced,
we deduce from (\ref{E:resol2}) the following inequality:
\begin{equation}\label{E:resol4}
(d-1)^2 - I \le t^1_X
\end{equation}
where:
\[
I = \begin{cases}\frac{(d-1)^2}{4}&\text{if $d$ is odd}\\ \frac{d(d-2)}{4} &\text{if $d$ is even}
\end{cases}
\]
These conditions easily imply the following result, whose part (1) is proved in a different way in \cite{aS05} and part (2) has been subsequently generalized in \cite{DS14} (see Remark \ref{R:free1} below).

\begin{prop}\label{P:free1}
(1) If $X$ is   nodal then it is not free unless $f=X_0X_1X_2$.

(2) If $X$ is irreducible, has $n$ nodes and $\kappa$ ordinary cusps as its only singularities and it is free 
then $\kappa \ge \frac{d^2}{4}$.
\end{prop}

\proof 1) If $X$ is nodal of degree $d=a+b+1$ then $t^1_X \le {a+b+1 \choose 2}$. It follows that
\[
(d-1)^2-t^1_X \ge (a+b)^2 - {a+b+1 \choose 2} = {a+b\choose 2} = ab+\frac{1}{2}[a(a-1)+b(b-1)]
\]
and this inequality is incompatible with the second condition (\ref{E:resol2}) unless
$a=b=1$.  This leaves  space for the existence of only one free 
(reducible) nodal curve:   the curve given by $f=X_0X_1X_2$, which is in fact free.   

2) Recalling that $t^1_X=n+2\kappa$ and combining  the inequality $n+\kappa \le {d-1\choose 2}$ with (\ref{E:resol4})    we obtain:
\[
(d-1)^2 - I \le \kappa + {d-1\choose 2}
\]
Now both possibilities for $I$ give the desired inequality after an easy calculation.
\qed

\begin{remark}\rm\label{R:free1}
In the recent preprint \cite{DS14} it has been proved that \emph{all} curves of degree $d \ge 4$ having only nodes and cusps are not free (see loc.cit., Example 4.5(ii)). The method of proof is
quite different, so we believe it can be useful to maintain the present weaker statement and its more elementary proof.

\end{remark} 

Several examples of free arrangements of lines are known. 
A notable example   is the dual of the configuration of flexes of a nonsingular plane cubic.
It consists of 9 lines meeting in 12 triple points.  Another free arrangement is given by $f=X_0X_1X_2(X_0-X_1)(X_1-X_2)(X_0-X_2)$: it has 4 triple points and 3 double points (see
\cite{sT12}, Ex. 3.4).  

The first example of free irreducible
plane curve has been given by Simis in \cite{aS05}. It is the sextic $X$ given by the polynomial:
\begin{equation}\label{Ex:simis}
f=4(X^2+Y^2+XZ)^3-27(X^2+Y^2)^2Z^2
\end{equation}
It has 4 distinct singular points, defined by the ideal
\[
\mathrm{rad}(J)= (YZ,2X^2+2Y^2-XZ) \]
One of them is a node 
and the other three are $E_6$-singularities. This curve is dual to a rational quartic $C$ with three nodes and three undulations (hyperflexes).
The $E_6$-singularities of $X$ are dual to the undulations of $C$. They have $\delta$-invariant $3$ and Tjurina number $6$.  Thus $t^1_X=3\cdot 6+1=19$.
Therefore $a+b=5$ and $ab=25-19=6$ and necessarily
\[
  T\langle X\rangle (-1) = \mathcal{O}(-3)\oplus \mathcal{O}(-2)
	\]
An	interesting example is the irreducible plane quintic  curve $X$ of equation $X_1^5-X_0^2X_2^3=0$. It has an $E_8$ and an $A_4$ singularity. They have respectively
$\delta = 4,2$ thus making the curve rational. On the other hand they have Milnor (equal to Tjurina) numbers equal to $8,4$ respectively, thus making
$t^1_X=12$.  The dual $X^\vee$ is again a quintic. According to (\ref{E:resol2}), if $X$ were free one should have  
\[
T\langle X\rangle(-1) = \mathcal{O}(-2)\oplus \mathcal{O}(-2)
\]
But $(f_0,f_1,f_2)$ has a linear syzygy (Example \ref{Ex:unstable}) and therefore this cannot be. 

 Other series of free irreducible plane curves are given in \cite{BC12,rN13,ST12,SiT12,sT12}. For a detailed discussion of freeness and more examples in the case of plane curves 
 we refer to \cite{DS14}.

\begin{example}\label{Ex:steinersur} \rm
The \emph{Steiner quartic surface} in $\P^3$, has equation in
normal (Weierstrass) form: $Z^2T^2+T^2Y^2+Y^2Z^2 = XYZT$. It is
irreducible and singular along the three coordinate axes for the
origin $(0,0,0,1)$. The jacobian ideal is
\[
J =(YZT,2YZ^2-XZT+2YT^2,2Y^2Z-XYT+2ZT^2,XYZ-2Y^2T-2Z^2T)
\]
and it turns out that $J^{sat}=J$. Therefore
$H^0_{\mathbf{m}}(R(f))=0$. Nevertheless it can be   computed  that $T\langle X \rangle$ is not free.
\end{example}

%%%%%%%%%%%%%%%%%%%%%%%%%%%%%%%%%%%%%%%%%%%%%%%%%%%

\noindent\textsc{address of the author:}

\noindent Dipartimento di Matematica e Fisica,
  Universit\`a Roma Tre \newline Largo S. L. Murialdo 1,
  00146 Roma, Italy.
  \smallskip
\newline  \texttt{sernesi@mat.uniroma3.it}

\end{document}